\documentstyle{article}

\newcommand{\sect}[1]{\section{#1}\setcounter{equation}{0}}

\newtheorem{th}{Theorem}[section]
\newtheorem{cor}[th]{Corollary}
\newtheorem{lem}[th]{Lemma}
\newtheorem{prop}[th]{Proposition}

\newtheorem{definition}[th]{Definition}

\def\sep{\;\vrule\;}
\def\proof#1. {\par
                      \ifdim\lastskip<15pt
                      \removelastskip\penalty-200
                      \vskip15pt plus3pt minus3pt
                      \fi
                       {\def\a{#1}
                       \ifx\a\empty
                       {\noindent\bf Proof.}
                       \else
                       {\noindent\bf Proof of #1.}
                       \fi}\enspace}
\def\restr#1{\,\vrule\,\lower1.75ex\hbox{$#1$}}

\def\endproof{\hfill\hspace{-6pt}\rule[-14pt]{6pt}{6pt}
\vskip22pt plus3pt minus 3pt}

\def\be{\begin{equation}}
\def\ee{\end{equation}}
\def\bea{\begin{eqnarray}}
\def\eea{\end{eqnarray}}
\def\bean{\begin{eqnarray*}}
\def\eean{\end{eqnarray*}}

\def\a{\alpha}
\def\b{\beta}
\def\d{\delta}

\def\e{\varepsilon}
\def\f{\varphi}
\def\F{\Phi}
\def\g{\gamma}

\def\i{\infty}
\def\k{\kappa}
\def\l{\lambda}
\def\L{\Lambda}

\def\s{\sigma}

\font\tenopen = cmbx10
\font\sevenopen = cmbx7
\font\fiveopen = cmbx5
\newfam\openfam
\def\open{\fam\openfam\tenopen}
\textfont\openfam = \tenopen
\scriptfont\openfam = \sevenopen
\scriptscriptfont\openfam = \fiveopen

\def\R{{\open R}}

\def\u{\underline}
\def\Q{{\bf Q}}
\def\Z{\bf Z}
\def\rz{{\R/\Z}}
\title{Weyl's uniform distribution under periodic perturbation\footnote{AMS Classification: 11K06; Keywords: Weyl's theory, density (mod 1),
uniform distribution (mod 1),
polynomial sequences, vector sequences, periodic perturbations}}
\author{Vilmos Totik}
\date{}
\begin{document}
\maketitle
\begin{abstract} We examine the uniform distribution theory of H. Weyl when
there is a periodic perturbation present. As opposed to the classical
setting, in this case the conditions for (mod 1) density and (mod 1)
uniform distribution turn out to be different.
\end{abstract}

\sect{Introduction}
In connection with asymptotics of orthogonal polynomials Peter Yuditskii \cite{Yud} asked if
the numbers $n\a-\cos(n\pi \a+\phi_0)$, $n=1,2,\ldots$, where $\a$ is irrational
 and $\phi_0$ is a given number, are dense modulo 1 (in short (mod 1)). Recall that
 if $\a$ irrational, then the sequence $n\a$, $n=1,2,\ldots$, is dense in
 $[0,1)$ (mod 1), so the problem asks if the same remains true if we perturb
 this sequence in the fashion described.

This paper has emerged from this question, and will examine  classical uniform distribution
theory when there is a periodic (or almost periodic) perturbation present.

We shall arrive at the most general results through simpler cases because those
simpler cases are interesting in themselves and because the proofs will be more
transparent  by not repeating simpler arguments as we proceed to higher generality.

Let us begin with a simple result with an elementary proof that answers Yuditskii's question.
\begin{th}\label{th1} If $f$ is a continuous periodic function and $\a$ is irrational,
then the sequence $\{n\a+f(n)\}_{n=1}^\i$ is dense in $[0,1)$ {\rm (mod 1)}.
\end{th}

Thus, mod 1  density is always the case, but we shall see later that (mod 1)  uniform
distribution may not be true. This is in contrast with the classical setting, where
the conditions for density and uniform distribution are mostly the same.

We shall have much more general density results, but for a start we give an elementary proof of
Theorem \ref{th1}.

\proof. Let $\e>0$ be small, and let $L$ be a large integer. Choose $\d<1/2$ such that $|f(u)-f(v)|<\e$ if
$|u-v|<\d$, and let $|f|\le M$, $M$ integer. If $\b$ is the/a period of $f$, then
choose a $k$ such that $k(\a/L)\in [-\e/L,\e/L]$ (mod 1)
and $k(1/\b)\in [-\d/\b,\d/\b]$ (mod 1)  (possible
by Dirichlet's simultaneous approximation theorem, see \cite[Theorem 201]{Hardy}). Then $k\a\in[-\e,\e]$ (mod $L$)
and $k\in [-\d,\d]$ (mod $\b$).

Assume, for example, that
$k\a\in [0,\e]$ (mod $L$), say $k\a=k_1L+\a_1$, where $k_1$ is an integer and
$\a_1\in (0,\e]$ (since $\a$ is irrational,
$\a_1$ cannot be 0). Consider the numbers
\[x_m=mk\a,\qquad 1\le m\le N,\]
where $N=[L/\a_1]$. We have
$x_m=m(k\a-k_1L)=m\a_1$ (mod $L$),
and these last numbers $m\a_1$ are equidistant on the interval $[0,L]$ with
neighboring distance $\a_1\le \e$. At the same time
\[|f(k(m+1))-f(km))|=|f(km+(k\  {\rm mod}\ \b))-f(km)|
\le \e\]
because of the choice of $\d$ and because of $(k \ {\rm mod}\ \b)\in [-\d,\d]$,
where $(k \ {\rm mod}\ \b)$ denotes the number $k'$ of smallest absolute
value for which $k=k'$ (mod $\b$).

Let $M_1=[M/\a_1]+1$. Since
\be M\le M_1\a_1\le M+\a_1<M+1\label{m1}\ee
 and
\be L-2-M<L-\a_1 -(M+\a_1)<(N-M_1)\a_1<L-M,\label{m2}\ee
 it follows that if
$y_m=x_m+f(km)$ (mod $L$), $y_m\in [0,L)$, then the sequence
$\{y_m\}_{m=M_1}^{N-M_1}$ (which is the same as
$\{m\a_1+f(km)\}_{m=M_1}^{N-M_1}$ in view of (\ref{m1}), (\ref{m2}) and $|f|\le M$) starts from a positive number $<2M+1$ and ends with a number
$>L-2M-2$, and the distance in between consecutive terms is at most $2\e$.
Therefore, for every $y\in [2M+1,L-2M-2]$ there is a term $y_m$ in the
sequence with $|y-y_m|\le 2\e$, and for $L>4M+5$ (which guarantees that
$[2M+1,L-2M-2]$ contains two consecutive integers) this shows that
to any point $y'\in [0,1)$ there is a $y_m$ closer than $2\e$ (mod 1), and (mod 1)  such a $y_m$ equals $km\a+f(km)$.\endproof

To illustrate that there is room for generalization, let us consider almost periodic functions.
 Bohr's almost periodic (real) functions (sometimes called uniformly almost periodic
 functions) coincide with the functions which can be uniformly approximated on $\R$ by generalized
 (real) trigonometric polynomials
 of the form $\sum c_je^{i\l_j x}$, see \cite[Section 84]{Bohr}\footnote{The original definition of almost periodicity is as follows:
 a continuous real function $f$ is said to be almost periodic if for every $\e>0$ there is an
 $L(\e)>0$ such for in every subinterval $I$ of $\R$ of length $\ge L(\e)$ there is a $t\in I$
 for which $|f(x+t)-f(x)|\le \e$ for all $x\in \R$. That the space of almost periodic functions
 coincides with the uniform closure of generalized trigonometric polynomials was one of Bohr's main theorems.}. Proposition \ref{th1} is the special case of

\begin{prop}\label{th2} If $f$ is almost periodic and $\a$ is irrational,
then  the sequence $\{n\a+f(n)\}_{n=1}^\i$ is dense in $[0,1)$ {\rm (mod 1)}.
\end{prop}

\proof. First of all notice that the proof of Proposition \ref{th1} easily implies the following:
if $f_1,\ldots,f_\nu$ are continuous periodic functions and $\a$ is irrational,
then the sequence $\{n\a+f_1(n)+\cdots+f_\nu(n)\}_{n=1}^\i$ is dense in $[0,1]$ (mod 1).
Indeed, let $\b_1,\ldots,\b_\nu$ be the periods of $f_1,\ldots,f_\nu$, respectively.
As before, let $\e>0$ be given, and let $L$ be again a large integer. Choose $\d$ such that $|f_\tau(u)-f_\tau(v)|<\e/\nu$
for all $\tau$ if
$|u-v|<\d$, and let $\sum_\tau|f_\tau|\le M$, $1\le \tau\le \nu$, $M$ integer.
Dirichlet's simultaneous approximation theorem can again be applied and we can choose
a $k$ such that $k(\a/L)\in [-\e/L,\e/L]$ (mod 1)
and $k(1/\b_\tau)\in [-\d/\b_\tau,\d/\b_\tau]$ (mod 1)  for all $\tau$.
The rest of the proof proceeds as before with minor modifications.

In particular, it follows that if $T=\sum_{\tau=1}^\nu (a_\tau\cos (\l_\tau x)+b_\tau
\sin (\l_\tau x))$ is a generalized trigonometric polynomial, then
the sequence $\{n\a+T(n)\}_{n=1}^\i$ is dense in $[0,1]$ (mod 1).
But every almost periodic function is the uniform
limit of such generalized trigonometric polynomials, and
the claim in the theorem follows: if $y\in (0,1)$ and $\e>0$ are given, then let
$\eta=\min\{\e,y,1-y\}$, choose
a $T$ with $|f-T|<\eta/2$ and then an $n$ with\footnote{$\{\cdot\}$ denotes fractional part.}
 $|y-\{n\a+T(n)\}|<\eta/2$, from which we obtain
 $|y-\{n\a+f(n)\}|<\eta<\e$.\endproof

\sect{Preliminaries, Weyl's theorems}\label{sectprel}
Let $\{\cdot\}$ denote fractional part. There are two slightly different ways
to talk about (mod 1)  values of a sequence $x_n$, $n=1,2,\ldots$, namely one can look
at the sequence $\{x_n\}$, $n=1,2,\ldots$, of the fractional parts as elements of $[0,1)$
or as elements of the torus $\R/{\bf Z}$ (i.e. when we identify $0$ and $1$ in $[0,1]$).
From the point of view of (mod 1)  denseness or (mod 1)  uniform distribution it does
not matter which view we take. In our discussion it will be more convenient to work
on $[0,1)$ and not on the torus because the proofs are easier that way and because
the perturbing function $f$ may not be conveniently reduced to a function on the torus.
We shall briefly discuss in Section \ref{secttorus}
one question when the torus vs. $[0,1)$ makes a difference.

Let $\{x_n\}_{n=1}^\i$ be a sequence of real numbers. This sequence is said to be uniformly
distributed on $[0,1)$ (mod 1) if for all $0\le a<b<1$ we have
\[\lim_{n\to\i}\frac{1}{n}\#\bigl\{{1\le k\le n\sep \{x_k\}\in [a,b]}\bigr\}=b-a.\]
This is the same that for every continuous function $f$ on $[0,1]$ we have
\[\lim_{n\to\i}\frac{1}{n}\sum_{k=1}^nf(\{x_k\})=\int_0^1 f.\]
In a similar fashion, if $\{X_n\}_{n=1}^\i$ is a sequence of vectors from $\R^d$, then
its (mod 1) uniform
distribution (on $[0,1)^d$) means that for all $0\le a_i<b_i<1$ we have
\[\lim_{n\to\i}\frac{1}{n}\#\bigl\{{1\le k\le n\sep \{X_k\}\in }\prod_{i=1}^d[a_i,b_i]\bigr\}=\prod_{i=1}^d(b_i-a_i),\]
which is the same that for every continuous function $f$ on $[0,1]^d$ we have
\[\lim_{n\to\i}\frac{1}{n}\sum_{k=1}^nf(\{x_k\})=\int_{[0,1]^d} f.\]

Weyl's first theorem \cite[Satz 1]{Weyl} (\cite[Theorem 2.1]{K}) on uniform
distribution asserts that $\{x_n\}_{n=1}^\i$  is uniformly distributed (mod 1)  on $[0,1)$ if and only if for all integers $a\not=0$
\[\lim_{n\to\i}\frac{1}{n}\sum_{k=1}^n e^{i2\pi a x_k}=0.\]
More generally, if
\be \{(x_{1,n},\ldots,x_{d,n})\}_{n=1}^\i\label{xnk}\ee
is a sequence from
$\R^d$, then
this sequence is uniformly distributed (mod 1) if and only if for all integer $d$-tuples
$(a_1,\ldots,a_d)$ with $\sum a_i^2\not=0$ we have
\[\lim_{n\to\i}\frac{1}{n}\sum_{k=1}^n e^{i2\pi (a_1x_{1,k}+\cdots+a_dx_{d,k})}=0,\]
see \cite[Satz 3]{Weyl} (\cite[Theorem 6.2]{K}).

In particular, (\ref{xnk}) is uniformly distributed on $[0,1)^d$ if an only if
every non-trivial integer combination sequence $y_n=a_1x_{1,n}+\cdots+a_dx_{d,n}$, $n=1,2,\ldots$,
is uniformly distributed on $[0,1)$.
This is not true if we  talk about density ((\ref{xnk}) may not be dense
in  $[0,1)^d$ even though every $\{y_n\}_{n=1}^\i$ is dense in $[0,1))$, but we shall see that for polynomial
sequences density and uniform distribution are the same.

Let $P_1,\ldots,P_d$ be real polynomials. It is clear that from the point of view
of the (mod 1)  density or (mod 1)  uniform distribution of the sequence
\[(P_1(n),\ldots,P_d(n)),\qquad n=1,2,\ldots,\]
 the constant terms in the polynomials are irrelevant, therefore
in what follows we shall always assume that the polynomials
are without constant term.

 For polynomial
sequences  Weyl proved in  \cite[Satz 9]{Weyl} (\cite[Theorem 3.2]{K}) that if $P$ (without constant term)
has at least one irrational coefficient, then $\{P(n)\}_{n=1}^\i$
is uniformly distributed on $[0,1)$ (mod 1). On the other hand, if all coefficients of $P$
are rational, then it is clear that $\{P(n)\}_{n=1}^\i$ has
only finitely many different terms (mod 1), so this sequence is
not dense in $[0,1)$ (mod 1).

To deal with the vector case we introduce the following definition.
\begin{definition} \label{def1} {\rm We call polynomials $P_1,\ldots,P_d$ (without constant term) $\Q$-independent, if
no non-trivial rational combination $r_1P_1+\cdots+r_d P_d$ is a polynomial
with rational coefficients. }
\end{definition}

Clearly, this is the same that no non-trivial
linear combination of $P_1,\ldots,P_d$ with integer coefficients
is a polynomial with integer coefficients.

Combining Weyl's theorems discussed so far we can see that  if $P_1,\ldots,P_d$ are
polynomials without constant term, then
\be \{(P_1(n),\ldots,P_d(n))\}_{n=1}^\i\label{xnkd}\ee
is  uniformly distributed in $[0,1)^d$ (mod 1)  precisely when $P_1,\ldots,P_d$
are $\Q$-inde\-pen\-dent, i.e.
   every non-trivial rational linear combination $\sum_j\a_j P_j$ has at least one
irrational coefficient. Furthermore, the condition is the same for the (mod 1)  density of
the sequence (\ref{xnkd}) in $[0,1)^d$. Indeed, if there is a non-trivial linear combination
$r_1P_1+\cdots+r_dP_d$ with integer coefficients which has integer coefficients, then
$r_1\{P_1(n)\}+\cdots+r_d\{P_d(n)\}$ is an integer $r$ lying in the interval
$[-(|r_1|+\cdots+|r_d|),|r_1|+\cdots+|r_d|]$. Therefore,  $(\{P_1(n)\},\ldots,\{P_d(n)\})$
lies in finitely many hyperplanes $r_1x_1+\cdots+r_dx_d=r$ of $\R^d$, hence this sequence
is not dense in $[0,1)^d$.

We can specialize these results to the case when $P_j(x)=\a_jx$ with some given $\a_j$ and  obtain that
the sequence
\be \{(n\a_1,\ldots,n\a_d )\}_{n=1}^\i\label{xnkv}\ee
of vectors is uniformly distributed (mod 1)  on $[0,1)^d$ if and only if
it is dense (mod 1)  in $[0,1)^d$, and this happens precisely
if the numbers $1,\a_1,\ldots,\a_d$ are rationally independent,
i.e. no non-trivial rational linear combination of them is 0.

\sect{Density}
In this section we shall deal with the problem discussed in the introduction
in more generality. We shall see that periodic or almost periodic perturbations
never prevent denseness (mod 1). The story will be different when we
shall deal with uniform distribution in the next section.

The following proposition is equivalent to Proposition \ref{th2}, but
 we present a formulation and a more streamlined proof for it that will be the
basis of the polynomial and vector cases to be followed.

\begin{prop} \label{th2*} If $\a$ is irrational and $f_1,\ldots,f_\nu$ are  almost periodic functions,
then  the sequence $\{n\a+\sum_jf_j(n)\}_{n=1}^\i$ is dense in $[0,1)$ {\rm (mod 1)}.
\end{prop}

\proof. As at the end of the proof of Proposition \ref{th2}, we may select each $f_j$
to be a trigonometric polynomial, and then, by considering individual terms
$\cos (\l x)$, $\sin (\l x)$ in the $f_j$'s and by increasing $\nu$
we may also assume that each $f_j$ is periodic
(say of the form $c\cos (\l x)$ or $c\sin (\l x)$ ). Let $\b_j$ be the period of
$f_j$, and set $\g_j=1/\b_j$. By replacing $f_j$ by $g_j(x)=f(\b_j x)$ it is sufficient
to show that if  each $g_j$ is periodic with period 1,
then  the sequence $\{n\a+\sum_jg_j(n\g_j)\}_{n=1}^\i$ is dense in $[0,1)$ (mod 1).

Let us choose a maximal subset $\{\g_\l\}_{\l\in \L}$ of $\{\g_1,\ldots,\g_\nu\}$
for which $\{1,\a\}\cup \{\g_\l\}_{\l\in \L}$ is a rationally independent set, i.e. 0 can be written
as a rational linear combination of these numbers only when all coefficients are 0. Thus, for every
$j\not\in \L$ there is a non-trivial dependence of the form
\[c_{j,0}\a+\sum_{\l\in \L} c_{j,\l}\g_\l+d_j\g_j=c_j,\]
where the $c$'s are integers and $d_j\not=0$ is an integer. With some large integer $L$
set $A=L\prod_{j\not \in \L} |d_j|$ (when $\L=\{1,\ldots,\nu\}$, then interpret the
product as 1), and consider numbers $n$ of the form
$n=Am$, $m=1,2,\ldots$. If $x\in (0,1)$ is given, then, by the rational independence of the set
$\{1,\a\}\cup \{\g_\l\}_{\l\in \L}$,
there is a subsequence ${\cal M}$ of the integers such that for $m\to\i$, $m\in {\cal M}$ we have
\[m\a\to x\ \  \mbox{(mod 1)}, \qquad m\g_\l \to \frac{1}{2}\ \  \mbox{(mod 1)}, \ \l\in \L.\]

Now
\[n\a+\sum_{j=1}^\nu g_j(n\g_j)=n\a+\sum_{j\in \L}g_j(n\g_j)+
\sum_{j\not \in \L}g_j\left( \frac{c_jn-c_{j,0}\a n-\sum_{\l\in \L} c_{j,\l}\g_\l n}{d_j}\right)\]
and this, by the 1 periodicity of the functions $g_j$, for $n=Am$ equals
\[A\{m\a\}+\sum_{j\in \L}g_j(A\{m\g_j\})+
\sum_{j\not \in \L}g_j\left(-c_{j,0}(A/d_j)\{m\a\}-\sum_{\l\in \L} c_{j,\l}(A/d_j)\{m\g_\l\}\right)\]
(mod 1), and this latter expression tends to
\be Ax+\sum_{j\in \L}g_j(A(1/2))+
\sum_{j\not \in \L}g_j\left(-c_{j,0}(A/d_j)x-\sum_{\l\in \L} c_{j,\l}(A/d_j)(1/2)\right)\label{89}\ee
as $m\to\i$, $m\in {\cal M}$.

From this we can conclude that the (mod 1) taken  sequence
\[\{n\a+\sum_{j=1}^\nu g_j(n\g_j)\}_{n=1}^\i\]
   is dense in the (mod 1)  taken range of the continuous function
\[h(x)=Ax+\sum_{j\in \L}g_j(A(1/2))+
\sum_{j\not \in \L}g_j\left(-c_{j,0}(A/d_j)x-\sum_{\l\in \L} c_{j,\l}(A/d_j)(1/2)\right)\]
defined on $(0,1)$, and it is sufficient to show that this (mod 1)  taken range is dense in $[0,1)$.

The function $h$ can also be considered on the interval $[0,1]$ where it is a continuous function, and,
 if $|f_j|\le M$, $1\le j\le \nu$, then
$h(0)\le \nu M$, while $h(1)\ge A-\nu M$. Therefore,  if $A>2\nu M+2$, then  its range contains
an interval $[B,B+1]$ for some integer $B$, and hence if we take its  range (mod 1), then we get the whole interval $[0,1)$.

This completes the proof.
\endproof

\subsection{The polynomial case}\label{pol}
As a generalization of what has been said before, we shall now examine the density of the
sequence
\be \Bigl\{P_0(n)+\sum_{j=1}^\nu f_j(P_j(n))\Bigr\}_{n=1}^\i,\label{00}\ee
(mod 1),  where $P_j$, $1\le j\le \nu$, are
polynomials. As before,  we may assume
that $P_j(0)=0$ for all $j$.

\begin{prop} \label{th2pol} Let $f_1,\ldots,f_\nu$ be almost periodic functions and
let $P_0,P_1,\ldots,P_\nu$ be polynomials without constant term. If $P_0$ has at least
one irrational coefficient, then the sequence {\rm (\ref{00})} is dense in $[0,1)$ {\rm (mod 1)}.
\end{prop}

\begin{cor} With the assumptions of the preceding proposition
\[\limsup_{n\to\i}\Bigl|\{P_0(n)\}-\bigl\{\sum_{j=1}^\nu f_j(P_j(n))\bigr\}\Bigr|\ge \frac{1}{2}.\]
\end{cor}

This corollary says that the fractional parts $\{P_0(n)\}$, $n=1,2,\ldots$, cannot be imitated by sums of the fractional
part of the values of almost periodic functions taken at polynomial places.
Note also that for $\nu=1$ and $f_1\equiv 1/2$ the equality is attained, so more than what is
claimed cannot be stated.

The corollary is immediate, for if the limsup in question was smaller than $1/2$, then
the fractional part of $P_0(n)-\sum_{j=1}^\nu f_j(P_j(n))$, $n=1,2,\ldots$,
would not accumulate at 1/2, so this sequence (mod 1)  would not be dense in $[0,1)$
contradicting Proposition \ref{th2pol}.

\proof Proposition \ref{th2pol}. We follow the preceding proof, and exactly as there we may assume that each $f_j$ is a continuous
periodic function with period 1.
Let us choose a maximal subset $\{P_\l\}_{\l\in \L}$ of $\{P_1,\ldots,P_\nu\}$
for which $\{P_0\}\cup \{P_\l\}_{\l\in \L}$ is a $\Q$-independent set. Thus, for every
$j\not\in \L$ there is a non-trivial dependence of the form
\[c_{j,0}P_0+\sum_{\l\in \L} c_{j,\l}P_\l+d_jP_j=R_j,\]
where the $c$'s are integers, $d_j\not=0$ is an integer and $R_j$ is a polynomial
with integer coefficients. With some large integer $L$
set again $A=L\prod_{j\not \in \L} |d_j|$ and define
\[P_j^*(x)=\frac{1}{A}P_j(Ax), \qquad R_j^*(x)=\frac{1}{A}R_j(Ax).\]
The $R_j^*$ has integer coefficients, and we claim that
$\{P_0^*\}\cup \{P_\l^*\}_{\l\in \L}$ is again a $\Q$-independent set.
Indeed, if this was not the case, then we could find integers $a_0, a_\l$, not all zero,
such that
\[Q(x)=a_0P_0^*(x)+\sum_{\l\in \L}a_\l P_\l^*(x)=
\frac{a_0}{A}P_0(Ax)+\sum_{\l\in \L}\frac{a_\l}{A} P_\l(Ax)\]
has integer coefficients. But then for $D$ equal to the degree of $Q$,
the coefficients of
\[a_0A^{D-1}P_0(x)+\sum_{\l\in \L}a_\l A^{D-1} P_\l(x)=A^D Q(x/A)\]
would have integer coefficients, which is impossible by the choice of
$\L$.

Consider numbers $n$ of the form
$n=Am$, $m=1,2,\ldots$. If $x\in (0,1)$ is given, then, by the $\Q$-independence of the set
$\{P_0^*\}\cup \{P_\l^*\}_{\l\in \L}$
and by Weyl's theorem,
there is a subsequence ${\cal M}$ of the integers such that for $m\to\i$, $m\in {\cal M}$, we have
\[P_0^*(m)\to x\ \  \mbox{(mod 1)}, \qquad P_\l^*(m) \to \frac{1}{2}\ \  \mbox{(mod 1)}, \ \l\in \L.\]
Now
\bean P_0(n)&+&\sum_{j=1}^\nu f_j(P_j(n))=P_0(n)+\sum_{j\in \L}f_j(P_j(n))\\
&+&
\sum_{j\not \in \L}f_j\left( \frac{R_j(n)-c_{j,0}P_0(n)-\sum_{\l\in \L} c_{j,\l}P_\l(n)}{d_j}\right).\eean
Here for all $j$ and $n=Am$ we have
\[P_j(n)=AP_j^*(m), \qquad R_j(n)=AR_j^*(m),\]
and since $A$ is divisible by $d_j$, it follows that each $R_j(n)/d_j$ has integer coefficients. Thus,
 by the 1-periodicity of the functions $f_j$, for $n=Am$ the sum
$ P_0(n)+\sum_{j=1}^\nu f_j(P_j(n))$  equals
\bean A\{P_0^*(m)\}&+&\sum_{j\in \L}f_j(A\{P_j^*(m)\})\\
&+&
\sum_{j\not \in \L}f_j\left(-c_{j,0}(A/d_j)\{P_0^*(m)\}-\sum_{\l\in \L} c_{j,\l}(A/d_j)\{P_\l^*(m)\}\right)\eean
(mod 1), and this latter expression tends to
\[Ax+\sum_{j\in \L}f_j(A(1/2))+
\sum_{j\not \in \L}f_j\left(-c_{j,0}(A/d_j)x-\sum_{\l\in \L} c_{j,\l}(A/d_j)(1/2)\right),\]
as $m\to\i$, $m\in {\cal M}$.

This is now the analogue of (\ref{89}), and from here the proof is completed exactly as
the preceding proof was completed from (\ref{89}).\endproof

Proposition \ref{th2pol} can be written in a seemingly more general, but actually equivalent form.

\begin{prop} \label{th2pol1} Let $F$ be an almost periodic function of $\nu$ variables\footnote{The set of these functions
$f: \R^\nu\to \R$
is the closure in the uniform norm of the algebra generated by the monomials $\cos \l x_j$, $\sin \l x_j$, $\l\in \R$, $1\le j\le \nu$.}
 and
let $P_0,P_1,\ldots,P_\nu$ be polynomials without constant term. If $P_0$ has at least
one irrational coefficient, then the sequence
\be \Bigl\{P_0(n)+F(P_1(n),\ldots,P_\nu(n))\Bigr\}_{n=1}^\i,\label{001}\ee
is dense in $[0,1)$ {\rm (mod 1)}.
\end{prop}

\proof. One can replace $F$ by a generalized trigonometric polynomial of several variables, and then one
can follow the preceding proof without much change.

An alternative is to use that if $F$ is a generalized trigonometric polynomial of $\nu$
variables, then
$F(P_1(n),\ldots,P_\nu(n))$ is of the form
\[\sum_{j=1}^{\nu_1} f_j(R_j(n))\]
where each $f_j$ is a $c\cdot\cos \l x$ or a $c\cdot \sin \l x$,
and where the $R_j$'s are polynomials (just use trigonometric
identities to convert products to sums).
This reduces then the claim to Proposition \ref{th2pol}.\endproof

\subsection{The vector case}
Let $\a_1,\ldots,\a_d$ be real numbers. Recall Weyl's theorem from Section \ref{sectprel}
according to which
for the (mod 1)  density in $[0,1)^d$ of the vectors
$(n\a_1,\ldots,n\a_d)$, $n=1,2,\ldots$, it is necessary and sufficient that  the
numbers $1,\a_1,\ldots,\a_d$
are rationally independent.

\begin{prop}\label{th1v} If $f_1,\ldots,f_d$ are almost periodic functions and
 the numbers $1,\a_1,\ldots,\a_d$ are rationally independent,
then the sequence
\be \bigl(n\a_1+f_1(n),\ldots,n\a_d+f_d(n)\bigr), \qquad n=1,2,\ldots,\label{seq}\ee
  is dense in $[0,1)^d$ {\rm (mod 1)}.
\end{prop}

\proof. The argument is familiar by now, see e.g. the proof of Proposition \ref{th2*}.

By replacing each $f_j$ by a trigonometric polynomial close to it, it is enough to prove that if
 $f_{i,j}$, $1\le i\le d$, $1\le j\le \nu$, are continuous periodic functions then for
any $\b_{i,j}$ the sequence
\be \Bigl(n\a_1+\sum_{j=1}^\nu f_{1,j}(n\b_{i,j}),\ldots,n\a_d+\sum_{j=1}^\nu f_{d,j}(n\b_{d,j})\Bigr), \qquad n=1,2,\ldots,\label{seq111}\ee
is dense in $[0,1)^d$ (mod 1). By applying dilation in $f_{i,j}$ and changing at the same
time the corresponding $\b_{i,j}$, we may assume that each $f_{i,j}$ is of period 1.

Select a maximal subset $\L\subset \{(i,j)\sep 1\le i\le d,\ 1\le j\le \nu\}$
such that
\be \{1,\a_1,\ldots,\a_d\}\cup \{\b_\l\sep \l\in \L\}\label{oi}\ee
is a rationally independent set. Then every $\b_{i,j}$ with $(i,j)\not\in \L$
can be expressed as a  rational combination
of the numbers in (\ref{oi}), say (with $\a_0=1$)
\[\b_{i,j}=\sum_{\tau=0}^d c_\tau^{i,j}\a_\tau+\sum_{\l\in \L}d_\l^{i,j} \b_\l.\]
If $q$ is the least common denominator of the coefficients in these representations,
then for some large $L$ set $A=Lq$ and consider the numbers $n=Am$, $m=1,2,\ldots$, for which
\[\Bigl(n\a_1+\sum_{j=1}^\nu f_{1,j}(n\b_{i,j}),\ldots,n\a_d+\sum_{j=1}^\nu f_{d,j}(n\b_{d,j})\Bigr)\]
is of the form
\be \Bigl(Am\a_1+\sum_{j=1}^\nu f_{1,j}(m(\cdots))\ldots,Am\a_d+\sum_{j=1}^\nu f_{d,j}(m(\cdots))\Bigr)\label{kk}\ee
where $(\cdots)$ represent integral combinations of the numbers in (\ref{oi}).

By Weyl's theorem, if $(x_1,\ldots,x_d)\in (0,1)^d$, then there is a subsequence ${\cal M}$ of the natural numbers such
that, as $m\to\i$, $m\in {\cal M}$, we have
\[\{m\a_j\}\to x_j,\quad  1\le j\le d, \qquad m\b_\l\to \frac12, \quad \l\in \L,\]
and we obtain from (\ref{kk}) that the sequence (\ref{seq}) is dense (mod 1)  in the (mod 1)  taken range of the
function
\be H(x_1,\ldots,x_d)=\Bigl(Ax_1+\sum_{j=1}^\nu f_{1,j}(\cdots)\ldots,Ax_d+\sum_{j=1}^\nu f_{d,j}(\cdots)\Bigr)\label{kk1}\ee
for $(x_1,\ldots,x_d)\in [0,1]^d$, where now $\cdots$ stand for some integral linear combinations of $x_1,\ldots,x_d, 1/2$.

Let $|f_{i,j}|\le M$ and $A>2\nu M+2$. By the following lemma (apply it to
\[h(x_1,\ldots,x_d)=H(x_1/A,\ldots,x_d/A)\]
and to the cube $Q=[0,A]^d$) the range of $H$ over
$[0,1]^d$ contains a cube of side-length 1, hence the (mod 1)  taken range of $H$ contains $[0,1)^d$,
from which the claim in the proposition follows.\endproof

In the next lemma we set $\underline x=(x_1,\ldots,x_d)$.
\begin{lem}\label{lem} If
\[h(\underline x)=(x_1+h_1(\underline x),\ldots,x_d+h_d(\underline x))\]
is a continuous mapping of a cube $Q\subset \R^d$ of side-length $A$ into $\R^d$, where
$|h_i|\le A/2-1$, then the range of $h$ contains a {\rm (}closed{\rm )} cube of side-length 1.
\end{lem}
\proof. This is an easy consequence of the Brouwer fixed point theorem. Indeed, we may assume
that $Q=[-A/2,A/2]^d$, and we show that then the range contains the cube $[-1/2,1/2]^d$.

Suppose to the contrary that this is not the case, and there is a point $\u y\in [-1/2,1/2]^d$
which is not in the range. For an $\u x\in Q$ connect
$h(\u x)$ with $\u y$, and let $\Phi(\u x)$ be the intersection of the extension
of this segment (beyond $\u y$) with the boundary $\partial Q$ of $Q$. This defines a continuous
mapping $\Phi: Q\to \partial Q$. Note that if $\u x$ belongs to the boundary of $Q$, say $x_1=-A/2$, then,
since $x_1+h_1(\u x)\le -1$, in forming $\Phi(\u x)$ we are connecting the point $h(\u x)$ lying in
the half-space $x_1\le -1$
with the point $\u y$ lying in the half-space $x_1\ge -1/2$. This  implies (say by
the convexity of the half-space $x_1<-1/2$) that the first coordinate of $\Phi(\u x)$ is
$\ge -1/2$, and hence $\Phi(\u x)\not= \u x$. Thus, $\F$ is a continuous mapping
of $Q$ into itself without a fixed point, which is impossible by the Brouwer fixed point theorem.

This contradiction proves the lemma.\endproof

\subsection{The polynomial-vector case}
The most general theorem in this paper on density is the following (cf. also Theorem \ref{thmaintor} below).
\begin{th}\label{thmain} Let $F_1,\ldots,F_d$ be
almost periodic functions of $\nu$ variables, and let
$P_{i,j}$ and $P_1,\ldots,P_d$ be polynomials without constant term.
If $\{P_1,\ldots,P_d\}$ are $\Q$-independent, then the sequence of vectors
\be \left\{\Bigl(P_1(n)+F_1(P_{1,1}(n)), \ldots,P_{1,\nu}(n)),\ldots,
P_d(n)+F_d(P_{d,1}(n),\ldots,P_{d\nu}(n))\Bigr)\right\}_{n=1}^\i\label{00101}\ee
is dense in $[0,1)^d$ {\rm (mod 1)}.
\end{th}

Recall that the $\Q$-independence of $\{P_1,\ldots,P_d\}$ is necessary
for density when each $F_i$ is identically 0.

\proof. All components of the proof have already been demonstrated in our previous discussions.

First of all, we may replace each $F_i$ by a trigonometric polynomial of $\nu$ variables (see the end of
the proof of Proposition \ref{th2}), then
using trigonometric identities, each $F_i(P_{i,1}(n),\ldots,P_{i,\nu}(n))$ can be written
in the form $\sum_{j=1}^\mu\f_i(S_{i,j})$ where each $\f_i$ is a $c\cdot \cos\l t$ or a  $c\cdot \sin\l t $
and $S_{i,j}$ are polynomials,
and finally such an expression is of the form   $\tilde F_i(S_{i,1}(n),\ldots,S_{i,\mu}(n))$,
where $\tilde F_i$ is a continuous function of $\mu$ variables which is periodic in each variable.
Thus, we may assume from the start that $F_i$ is periodic in each of its variable, and by scaling
we may even assume that the periods are 1.

Select a maximal subset $\L\subset \{(i,j)\sep 1\le i\le d,\ 1\le j\le \nu\}$
such that
\be \{P_1,\ldots,P_d\}\cup \{P_\l\sep \l\in \L\}\label{oi1}\ee
is a $\Q$-independent set of polynomials. Then for every $P_{i,j}$ with $(i,j)\not\in \L$
there is a relation with the polynomials  in (\ref{oi1}) of the form
\[\sum_{\tau=1}^d c_\tau^{i,j}P_\tau+\sum_{\l\in \L} c_\l^{i,j}P_\l+d_{i,j}P_{i,j}=R_{i,j},\]
where the $c$'s are integers, $d_{i,j}\not=0$ is an integer and $R_{i,j}$
is a polynomial with integer coefficients. We consider again
with some large integer $L$ the number
$A=L\prod_{(i,j)\not \in \L} |d_{i,j}|$ and define for a polynomial $P$ the polynomial
$P^*(x)=\frac{1}{A}P(Ax)$. The $*$-transform of the system (\ref{oi1}) is again
$\Q$-independent (see the proof of Proposition \ref{th2pol}), so, by Weyl's theorem,
if $(x_1,\ldots,x_d)\in (0,1)^d$ is given, then
there is a subsequence ${\cal M}$ of the integers such that for $m\to\i$, $m\in {\cal M}$ we have
\[P_j^*(m)\to x_j\  \mbox{(mod 1)}, \  j=1,\ldots,d,\qquad P_\l^*(m) \to
\frac{1}{2}\ \  \mbox{(mod 1)}, \ \l\in \L.\]
By considering numbers of the form $n=Am$, $m=1,2,\ldots$, and follow the steps in the proof of
Propositions \ref{th2pol} and \ref{th1v} we can conclude that the closure of the (mod 1)
taken sequence (\ref{00101}) includes the (mod 1)  taken range of the function
\[ H(x_1,\ldots,x_d)=\Bigl(Ax_1+F_1(\cdots)\ldots,Ax_d+ F_d(\cdots)\Bigr)\]
for $(x_1,\ldots,x_d)\in [0,1]^d$, where, for each of the variables
of each $F_i$, the $\cdots$ stands for some integral
linear combinations of $x_1,\ldots,x_d, 1/2$. That this (mod 1)  taken range is
$[0,1)^d$ follows from Lemma \ref{lem} precisely as in the proof of Proposition \ref{th1v}.\endproof

\sect{Uniform distribution}
In this section we shall examine
the uniform distribution of the sequences that have been considered before.
In the classical case Weyl's theorems implies that the conditions for
denseness and uniform distribution are the same. We shall see that
when a periodic perturbation is present the situation is
different: denseness is always true (as we have seen in the preceding
section), but that is not the case for uniform distribution.

\subsection{The scalar case} Let us start with the simplest possible situation.
\begin{prop} \label{th3} Let $\a$ be irrational and let $f$ be a $\b$-periodic continuous function. Then
the sequence $\{n\a+f(n)\}_{n=1}^\i$ is uniformly distributed {\rm (mod 1)}  if  $\b$ is rational
or $1,\a,1/\b$ are rationally independent.
\end{prop}

Let us mention that Proposition \ref{th3} is precise in the sense that if $\b$ is irrational and
$1,\a,1/\b$ are rationally dependent, then there is a $\b$-periodic continuous function $f$ such that
the sequence $\{n\a+f(n)\}_{n=1}^\i$ is not uniformly distributed (mod 1).
See Proposition \ref{th4}.

\proof. Set  $\g=1/\b$ and $g(x)=f(\b x)$. Then $g$ is 1-periodic
 and $n\a+f(n)=\{n\a\}+g(\{n\g\})$ (mod 1).

Consider first the case when $\g=1/\b$ is rational, say $\g=p/q$. Then for $n=mq+s$, $0\le s<q$, we have
\[n\a+g(n\g)=\{m(q\a)\}+(g(s\g)+s\a)\qquad \mbox{(mod 1)},\]
 and since here the second term on the right is independent of
$m$ and $q\a$ is irrational, it follows from the uniform distribution of $\{m(q\a)\}_{m=1}^\i$
(mod 1)  that
the sequence $\{(mq+s)\a+f((mq+s)\a)\}_{m=1}^\i$ is uniformly distributed (mod 1). Since this is true
for all $0\le s<q$, the (mod 1)  uniform distribution of $\{n\a+f(n)\}_{n=1}^\i$ follows.

Let us now consider the case when $1,\a,1/\b$, i.e. $1,\a,\g$ are rationally independent. Then
the vector sequence $\{(n\a,n\g)\}_{n=1}^\i$ is uniformly distributed
in $[0,1)^2$ (mod 1).

Let $0<a<b<1$ be fixed, and for an $0<\e<\min\{a,1-b\}$ let $m$ be so large that the
oscillation\footnote{Here and in what follows, the oscillation of a function on a set
is the difference of its supremum and infimum on that set} of $g$ over any interval of length $<1/m$ is smaller
than $\e$.  For $s=0,1,\ldots,m-1$ and for an $n$ consider the sets
\[K_{1,s}=\left\{1\le k\le n\sep k\g\in \left[\frac{s}{m},\frac{s+1}{m}\right],\
k\a\in \left[a-g\left(\frac{s}{m}\right)+\e,b-g\left(\frac{s}{m}\right)-\e\right]\right\},\]
and
\[K_{2,s}=\left\{1\le k\le n\sep k\g\in \left[\frac{s}{m},\frac{s+1}{m}\right],\
k\a\in \left[a-g\left(\frac{s}{m}\right)-\e,b-g\left(\frac{s}{m}\right)+\e\right]\right\},\]
where the relations are understood (mod 1). Since for
\[k\g\in \left[\frac{s}{m},\frac{s+1}{m}\right]\qquad \mbox{(mod 1)}\]
we have
\[g(k\g)\in \left[g\left(\frac{s}{m}\right)-\e,g\left(\frac{s}{m}\right)+\e\right]\]
by the choice of $m$, it follows that for $k\in K_{1,s}$ we have
$k\a+g(k\g)\in [a,b]$ (mod 1), and conversely, if
$k\g\in \left[\frac{s}{m},\frac{s+1}{m}\right]$ and
$k\a+g(k\g)\in [a,b]$ (mod 1), then
$k\in K_{2,s}$. By the (mod 1)  uniform distribution of the sequence of vectors $\{(k\a,k\g)\}_{k=1}^\i$ we have
\[\#K_{1,s}=(1+o(1))n\frac{b-a-2\e}{m},\qquad \#K_{2,s}=(1+o(1))n\frac{b-a+2\e}{m},\]
hence it follows by summing these for all $s=0,1,\ldots,m-1$ that
\[\liminf_{n\to\i}\frac1{n}\#\left\{1\le k\le n\sep k\a+g(k\g)\in [a,b] \mbox{\ \ (mod 1)}\right\}\ge b-a-2\e\]
and
\[\limsup_{n\to\i}\frac1{n}\#\left\{1\le k\le n\sep k\a+g(k\g)\in [a,b] \mbox{\ \ (mod 1)}\right\}\le b-a+2\e.\]
Since $\e>0$ and $0<a<b<1$ are arbitrary, the uniform distribution of $\{n\a+f(n)\}_{n=1}^\i$ (mod 1)  follows.\endproof

We shall need the extension of Proposition \ref{th3} to more than one term. To this end, note that
for an irrational $\a$ the condition
``$\b$ is rational or $1,\a,1/\b$ are rationally independent" can be phrased as
``$\a$ cannot be written as a rational linear combination of $1$ and $1/\b$". It turns out that this is the right
formulation in the general case.

\begin{prop} \label{th31} Let $\a$ be irrational and let $f_1,\ldots,f_\nu$ be continuous periodic functions with periods
$\b_1,\ldots,\b_\nu$. If $\a$ cannot be written as a rational linear combination of the numbers
$1,1/\b_1,\ldots,1/\b_\nu$, then
the sequence $\{n\a+f_1(n)+\cdots f_\nu(n)\}_{n=1}^\i$ is uniformly distributed {\rm (mod 1)}.
\end{prop}

This is best possible in the sense that if $\a$ and $\b_j$ are given and
$\a$ can be written as a
rational linear combination of the numbers
$1,1/\b_1,\ldots,1/\b_\nu$, then there are $\b_j$-periodic continuous functions $f_j$ for $j=1,\ldots,\nu$,
such that the sequence $n\a+f_1(n)+\cdots f_\nu(n)$, $n=1,2,\ldots$,
is not uniformly distributed (mod 1).
See Proposition \ref{th4}.

\proof. A large part of the proof has been already established in the preceding proof.

Set  $\g_i=1/\b_i$ and $g_i(x)=f_i(\b_i x)$. Then $g_i$ are 1-periodic
 and
 \[n\a+f_1(n)+\cdots f_\nu(n)=\{n\a+g_1(n\g_1)+\cdots g_\nu(n\g_\nu)\}\qquad \mbox{(mod 1)}.\]

Select a maximal subset of $\g_1,\ldots,\g_\nu$ that forms with 1 a rationally independent
set. We may assume that $\{\g_1,\ldots,\g_\mu\}$ is this maximal subset. Then, by the assumption of
the theorem, $1,\a,\g_1,\ldots,\g_\mu$ are rationally independent, furthermore every $\g_j$ with $j>\mu$ is a rational
combination of the form
\be \g_j=d_j+\sum_{i=1}^\mu c_{j,i}\g_i.\label{tt00}\ee
Let $q$ be a common multiple of the denominators of these $d_j$. If $k$ is of the form
$k=ql+r$, $0\le r<q$, then
\[k\g_j=rd_j+\sum_{i=1}^\mu c_{j,i}k\g_i\quad \mbox{(mod 1)},\]
and it is sufficient to show the uniform distribution in question along
every sequence ${\cal N}_r=\{ql+r\}_{l=1}^\i$ where $0\le r<q$ is fixed.

For a large $m$ partition $[0,1]^\mu$ by the hyperplanes $x_i=t/m$, $t=1,\ldots,m-1$, $1\le i\le \mu$,
 in $\R^\mu$ into $m^\mu$
little cubes $I_s$, $1\le s\le m^\mu$. The image of such an $I_s$ under the mapping
\be (x_1,\ldots,x_\mu)\to rd_j+\sum_{i=1}^\mu c_{j,i}x_i\label{tt100}\ee
(cf. (\ref{tt00})) is an interval $I_{j,s}$ of length $\le C/m$.

Now let  $0<a<b<1$ be fixed, and for an $0<\e<\min\{a,1-b\}$ let $m$ be so large that the oscillation of
$g_1(x_1)+\cdots + g_\mu(x_\mu)$  over any $I_s$, as well as the oscillation of each
$g_j$, $j>\mu$, over any $I_{j,s}$ is smaller than $\e/2\nu$. Let, furthermore, $Q_s=(Q_{1.s},\ldots,Q_{\mu,s})$ be an
arbitrary point from $I_s$, let $\tilde Q_{j,s}$ be the image of $Q_s$ under the mapping (\ref{tt100}),
and set
\[T_s=\sum_{i=1}^\mu g_i(\g_i Q_{i,s})+\sum_{j=\nu+1}^d  g_j(\g_j \tilde Q_{j,s}).\]

 For $s=1,\ldots,m^\mu$ consider the sets
\[K_{1,s}=\left\{1\le k\le n,\ k\in {\cal N}_r\sep (k\g_1,\ldots,k\g_\mu)\in I_s,\
k\a\in [a-T_s+\e,b-T_s-\e]\right\},\]
and
\[K_{2,s}=\left\{1\le k\le n,\ k\in {\cal N}_r\sep (k\g_1,\ldots,k\g_\mu)\in I_s,\
k\a\in [a-T_s-\e,b-T_s+\e]\right\},\]
where the relations are understood (mod 1). Exactly as in the preceding proof,
by the choice of $m$, it follows that for $k\in K_{1,s}$ we have
\be k\a+g_1(\g_1 k)+\cdots g_\nu(\g_\nu k)\in [a,b] \quad {\rm (mod}\  1),\label{fgh}\ee
 and conversely, if
$(k\g_1,\ldots,k\g_\mu)\in I_s$ and (\ref{fgh}) is true, then
$k\in K_{2,s}$.
By the (mod 1)  uniform distribution of the sequence of vectors
\[\{(k\a,k\g_1,\ldots,k\g_\mu)\}_{k=1,\; k\in{\cal N}_r}^\i\]
(which follows from the rational independence of $(1,\a,\g_1,\ldots,\g_\mu)$), we have
\[\#K_{1,s}=(1+o(1)n\frac{b-a-2\e}{q\cdot m^\mu},\qquad \#K_{2,s}=(1+o(1))n\frac{b-a+2\e}{q\cdot m^\mu}.\]
It follows by summing these for all $s=1,\ldots,m^\mu$ that
\bean \liminf_{n\to\i}\frac{q}{n}\#\left\{1\le k\le n,\ k\in {\cal N}_r\sep k\a+g_1(\g_1 k)+\cdots g_\nu(\g_\nu k)\in [a,b]
\mbox{\ \ (mod 1)}\right\}\\
\ge b-a-2\e\eean
and
\bean \limsup_{n\to\i}\frac{q}{n}\#\left\{1\le k\le n,\  k\in {\cal N}_r\sep k\a+g_1(\g_1 k)+\cdots g_\nu(\g_\nu k)\in [a,b]
\mbox{\ \ (mod 1)}\right\}\\
\le b-a+2\e.\eean
Since $\e>0$ and $0<a<b<1$ are arbitrary and since there are $n/q+O(1)$ numbers $k\in {\cal N}_r$ with
$1\le k\le n$, the uniform distribution of
\[\{k\a+g_1(\g_1 k)+\cdots g_\nu(\g_\nu k)\}_{k=1,\;k\in {\cal N}_r}^\i\]
 (mod 1)  follows.\endproof

Next we state a somewhat more general form (but actually, since the number $\nu$ of terms in Proposition
\ref{th31} is arbitrary, and since, by Weierstrass theorem, every periodic continuous function
of several variables can be uniformly approximated by trigonometric polynomials of the same variables,
one can show that this seemingly more general form is equivalent to
Proposition \ref{th31}).

\begin{prop} \label{th31000} Let $\a$ be irrational and let $F$ be a continuous
function of the variables $x_1,\ldots,x_\nu$ which is periodic in each variable with periods
$\b_1,\ldots,\b_\nu$, respectively. If $\a$ cannot be written as a rational linear combination of the numbers
$1,1/\b_1,\ldots,1/\b_\nu$, then
the sequence $\{n\a+F(n,\ldots,n)\}_{n=1}^\i$ is uniformly distributed {\rm (mod 1)}.
\end{prop}

In particular, if
\[F(x_1,\ldots,x_\nu)=f_1(x_1)+\cdots+f_\nu(x_\nu),\] where each $f_j$ is a continuous periodic function
of period $\b_j$, then we obtain Proposition \ref{th31}.

\proof. We follow the preceding proof.

Set  $\g_i=1/\b_i$ and
\[G(x_1,\ldots,x_\nu)=F(\b_1 x_1,\ldots,\b_\nu x_\nu).\]
 Then $G$ is 1-periodic in each variable
 and
 \[n\a+F(n,\ldots,n)=\{n\a+G(\{n\g_1\},\ldots,\{n\g_\nu\})\}\qquad \mbox{(mod 1)}.\]

Select again a maximal subset of $\g_1,\ldots,\g_\nu$ that forms with 1 a rationally independent
set. We may assume that $\{\g_1,\ldots,\g_\mu\}$ is this maximal subset. Then, by the assumption of
the theorem, $1,\a,\g_1,\ldots,\g_\mu$ are rationally independent, furthermore every $\g_j$ with $j>\mu$ is a rational
combination of the form (\ref{tt00}).
For a large $m$ partition $[0,1]^\mu$ by the hyperplanes $x_i=t/m$, $t=1,\ldots,m-1$ in $\R^\mu$ into $m^\mu$
little cubes $I_s$, $1\le s\le m^\mu$, as before. The image of such an $I_s$ under the mapping (\ref{tt100})
(cf. (\ref{tt00})) is again an interval $I_{j,s}$ of length $\le C/m$.

Now  for fixed $0<a<b<1$ and for an $0<\e<\min\{a,1-b\}$ let $m$ be so large that the oscillation of
$G$ over any set of diameter $\le C\nu/m$ is smaller than $\e/2$. Let furthermore, $Q_s=(Q_{s,1},\ldots,Q_{s,\mu})$  be an
arbitrary point from $I_s$, let $\tilde Q_{j,s}$ be the image of $Q_s$ under the mapping (\ref{tt100}),
and set
\[T_s=G(\g_1 Q_{s,1},\ldots,\g_\mu Q_{s,\mu},\g_{\mu+1}\tilde Q_{\mu+1,s},\ldots, \g_{\nu}\tilde Q_{\nu,s}).\]
With these modifications for the sets
\[K_{1,s}=\left\{1\le k\le n,\ k\in{\cal N}_r\sep (k\g_1,\ldots,k\g_\mu)\in I_s,\
k\a\in [a-T_s+\e,b-T_s-\e]\right\},\]
and
\[K_{2,s}=\left\{1\le k\le n,\ k\in{\cal N}_r\sep (k\g_1,\ldots,k\g_\mu)\in I_s,\
k\a\in [a-T_s-\e,b-T_s+\e]\right\},\]
(where the relations are understood (mod 1)) we get again  that for $k\in K_{1,s}$ we have
$k\a+G(\g_1 k,\cdots ,\g_\nu k)\in [a,b]$ (mod 1), and conversely, if
$(k\g_1,\ldots,k\g_\mu)\in I_s$ and
$k\a+G(\g_1 k,\cdots,\g_\nu k)\in [a,b]$ (mod 1), then
$k\in K_{2,s}$.
By the (mod 1)  uniform distribution of the sequence of vectors $\{(k\a,k\g_1,\ldots,k\g_\mu)\}_{k=1,\;k\in {\cal N}_r}^\i$
(which is a consequence of the  rational independence of $(1,\a,\g_1,\ldots,\g_\mu)$), we get again
\[\#K_{1,s}=(1+o(1)n\frac{b-a-2\e}{q\cdot m^\mu},\qquad \#K_{2,s}=(1+o(1))n\frac{b-a+2\e}{q\cdot m^\mu},\]
and from here the uniform distribution of
\[\{k\a+G(k\g_1,\cdots ,k\g_\nu)\}_{k=1,\;k\in{\cal N}_r}^\i\]
 (mod 1) for each $0\le r<q$ follows as before, and that proves
Proposition \ref{th31000}.\endproof

\subsection{The polynomial case}
Now we extend Proposition \ref{th31000} to the polynomial case. As before, we may
assume without loss of generality that the polynomials we are dealing with are without constant term.

\begin{prop} \label{th312} Let $P_0,P_1,\ldots,P_d$ be polynomials without constant term and
let $F$ be a continuous
function of the variables $x_1,\ldots,x_\nu$ which is periodic in each variable with periods
$\b_1,\ldots,\b_\nu$, respectively.  If $P_0$ cannot be written as a rational linear combination of the polynomials
$P_1/\b_1,\ldots,P_d/\b_d$ and of a polynomial with rational coefficients, then
the sequence
\be \{P_0(n)+F(P_1(n),\cdots,P_\nu(n))\}_{n=1}^\i\label{subs}\ee
 is uniformly distributed {\rm (mod 1)}.
\end{prop}
Note that Proposition \ref{th31000} is the special case of this when $P_0(x)=\a x$ and
$P_j(x)=x$ for all other $j$.

In particular, if
\[F(x_1,\ldots,x_\nu)=f_1(x_1)+\cdots+f_\nu(x_\nu),\]
where each $f_j$ is a continuous periodic function
of period $\b_j$, then under the conditions of the theorem the (mod 1)  density of the sequence
\[\{P_0(n)+f_1(P_1(n))+\cdots+f_\nu(P_\nu(n))\}_{n=1}^\i\]
 follows.
In Proposition \ref{th41} we shall show that this is best possible in the sense that if $P_0,P_1,\ldots,P_\nu$
are given and $P_0$ can be written as a rational linear combination of the polynomials
$P_1/\b_1,\ldots,P_\nu/\b_\nu$ and of a polynomial with rational coefficients, then there are
continuous functions $f_j$ of period $\b_j$, $1\le j\le \nu$, such that
the sequence
\[P_0(n)+f_1(P_1(n))+\cdots f_d(P_\nu(n)), \qquad n=1,2,\ldots,\]
 is not
uniformly distributed (mod 1).

\proof. The proof is almost identical to that of Proposition \ref{th31}. Indeed, as there
set  $\g_i=1/\b_i$ and
\[G(x_1,\ldots,x_\nu)=F(\b_1 x_1,\ldots,\b_\nu x_\nu).\]
Then  $F(P_1(n),\ldots,P_\nu(n))=G(\g_1 P_1(n),\ldots,\g_\nu P_\nu(n))$.
Select a maximal subset of $\g_1P_1,\ldots,\g_\nu P_\nu$  for which no non-trivial linear combination with integer coefficients produces
a polynomial with integer coefficients. We may assume that $\g_1P_1,\ldots,\g_\mu P_\mu$ is this maximal subset.
Then, by the assumption of
the theorem, $P_0,\g_1P_1,\ldots,\g_\mu P_\mu$ is $\Q$-independent (see Definition \ref{def1}), but every $\g_jP_j$ with $j>\nu$ is a
combination of the form
\be \g_jP_j(x)=R_j(x)+\sum_{i=1}^\mu c_{j,i}\g_iP_i(x),\label{tt11}\ee
where the $c_{j,i}$ are rationals and $R_j$ has rational coefficients.
Let $q$ be a common multiple of the denominators of all the coefficients
in all $R_j$, $\mu< j\le \nu$. It is sufficient to verify the uniform
distribution of each subsequence of (\ref{subs}) for which the indices $n$
belong to ${\cal N}_r=\{ql+r\}_{l=1}^\i$ with a fixed $0\le r<q$. If
$k=ql+r$ is such an index, then
\[ \g_jP_j(k)=R_j(r)+\sum_{i=1}^\mu c_{j,i}\g_iP_i(k)\qquad \mbox{(mod 1)},\]
and on the right $R_j(r)$ is independent of $l=1,2,\ldots$. This $R_j(r)$
plays the role of $rd_j$ in (\ref{tt100}).
Since
$P_0,\g_1P_1,\ldots,\g_\mu P_\mu$ are $\Q$-independent, it follows that the polynomials
$P_0(qx+r),\g_1P_1(qx+r)\ldots,\g_\mu P_\mu(qx+r)$ are also $Q$-independent (note that if a polynomial
$R$ has rational coefficients then so does $R((y-r)/q)$). Hence Weyl's theorem ensures that
the sequence
\[(P_0(k),\g_1P_1(k),\ldots,\g_\mu P_\mu(k)),\qquad k=ql+r,\ l=1,2,\ldots,\]
of vectors is uniformly distributed (mod 1)  on $[0,1)^{\nu+1}$.
This is the analogue of the uniform distribution of the sequence
$\{(k\a,n\g_1,\ldots,k\g_\mu)\}_{k=1,\;k\in{\cal N}_r}^\i$ in the preceding proof
(while (\ref{tt11}) is the analogue of (\ref{tt00})), and by replacing the
latter by the previous one, that proof goes through without much change. \endproof

\subsection{The vector case}
Now we shall consider the vector case, namely the uniform distribution
of vector sequences of the form
\be \Bigl\{\bigl(n\a_1+f_1(n), \ldots,n\a_d+f_d(n)\bigr)\Bigr\}_{n=1}^\i\label{vvec}\ee
 (mod 1), where the $f_i$ are continuous periodic functions. To this end we introduce the following
 definition.

 \begin{definition} {\rm We say that a system $(\a_1,\b_1),\ldots,(\a_d,\b_d)$ of number pairs is totally
 $\Q$-independent, if for any non-empty subset $\L\subset \{1,\ldots,d\}$
 it is true, that if $a_\l$, $\l\in \L$, are {\it nonzero} integers, then
 the number $\sum_{\l\in \L} a_\l \a_\l$ cannot be written as a rational linear combination of the numbers
 $1,\b_\l$, $\l\in \L$.}
 \end{definition}

As an example consider the pairs
$(\sqrt 2,\sqrt 3)$, $(\sqrt 3,\pi)$, $(\pi^2,\sqrt 2)$. This is easily seen
to be totally $\Q$-independent, though $\a_1+\a_2=\sqrt 2+\sqrt 3=\b_1+\b_3$ (note that total $\Q$-independence requires
that $\a_1+\a_2=\sqrt 2+\sqrt 3$ should not be a rational
linear combination of $\b_1=\sqrt 2$ and $\b_2=\pi$ -- we cannot
use in this combination $\b_3$).

\begin{prop} \label{th311} Let $f_1,\ldots,f_d$ be continuous periodic functions with periods
$\b_1,\ldots,\b_d$. If $(\a_1,1/\b_1),\ldots,(\a_d,1/\b_d)$ are totally $\Q$-independent,
 then
the sequence {\rm (\ref{vvec})} is uniformly distributed on $[0,1)^d$ {\rm (mod 1)}.
\end{prop}

In Proposition \ref{Th42} we shall show that one cannot do more than Proposition \ref{th311},
for  uniform distribution is not true
(for some functions $f_j$ of period $\b_j$) if
$(\a_1,1/\b_1),\ldots,(\a_d,1/\b_d)$ are not totally $\Q$-independent.

\proof. By Weyl's criterion, a sequence
\be (v_{1,n},\ldots,v_{d,n}),\qquad n=1,2,\ldots,\label{df}\ee
of vectors is uniformly distributed on $[0,1)^d$ (mod 1)  precisely if for all integers $a_1,\ldots,a_d$, $\sum a_j^2>0$,
the sequence of numbers
\be a_1v_{1,n}+\cdots+a_nv_{d,n},\qquad n=1,2,\ldots,\label{df1}\ee
is uniformly distributed on $[0,1)$ (mod 1). Apply this to
$v_{j,n}=\a_j n+f_j(n)$. If in the combination
the non-zero $a_j$ are $a_\l$, $\l\in \L$, then we need to prove the (mod 1)  uniform
distribution on $[0,1)$ of the sequence
\[\left(\sum_{\l\in \L} a_\l \a_\l\right) n+\sum_{\l\in \L} a_\l f_\l(n),\qquad n=1,2,\ldots,\]
which follows from Proposition \ref{th31} because, by total $\Q$-independence,
$\sum_{\l\in \L} a_\l \a_\l$ cannot be written as a rational linear combination of the numbers
 $1, 1/\b_\l$, $\l\in \L$.\endproof

We can state a more general form of the preceding proposition. To this end  we introduce the following
 definition.

 \begin{definition} {\rm We say that a system
 \[(\a_1,\b_{1,1},\ldots,\b_{1,\nu}),\ldots,(\a_d,\b_{d,1},\ldots,\b_{d,\nu})\]
  of $(\nu+1)$-tuples of numbers is totally
 $\Q$-independent, if for any non-empty subset $\L\subset \{1,\ldots,d\}$
 it is true, that if $a_\l$, $\l\in \L$, are {\it nonzero} integers, then
 the number $\sum_{\l\in \L} a_\l \a_\l$ cannot be written as a rational linear combination of the numbers
 $1,\b_{\l,j}$, $\l\in \L,1\le j\le \nu$.}
 \end{definition}

\begin{prop} \label{th311*} Let $F_1,\ldots,F_d$ be continuous periodic functions of $\nu$ variables
such that $F_i$ is periodic in the $j$-th variable with period
$\b_{i,j}$. If
 \[(\a_1,1/\b_{1,1},\ldots,1/\b_{1,\nu}),\ldots,(\a_d,1/\b_{d,1},\ldots,1/\b_{d,\nu})\]
 are totally $\Q$-independent,  then
the sequence
\be \Bigl\{\bigl(n\a_1+F_1(n,\ldots,n), \ldots,n\a_d+F_d(n,\ldots,n)\bigr)\Bigr\}_{n=1}^\i\label{vvecv}\ee
of vectors is uniformly distributed on $[0,1)^d$ {\rm (mod 1)}.
\end{prop}

\proof. Just follow the proof of Proposition \ref{th311}. In this case we apply
Weyl's criterion to a linear combination $\sum_j a_j v_{j,n}$ with
$v_{j,n}=\a_j n+F_j(n,\ldots,n)$. If in the combination
the non-zero $a_j$ are $a_\l$, $\l\in \L$, then we need to prove the (mod 1)  uniform
distribution on $[0,1)$ of the sequence
\[\left(\sum_{\l\in \L} a_\l \a_\l\right) n+\sum_{\l\in \L} a_\l F_\l(n,\ldots,n),\qquad n=1,2,\ldots,\]
which follows from  Proposition \ref{th31000} because, by total $\Q$-independence,
$\sum_{\l\in \L} a_\l \a_\l$ cannot be written as a rational linear combination of the numbers
$1,\  1/\b_{\l,j}, \l\in \L,\ 1\le j\le \nu$. Indeed, if, say, $\L=\{1,\ldots,\k\}$ (which
we may assume), then one should replace $\nu$ in Proposition \ref{th31000}  by $\k\nu$ and
apply Proposition \ref{th31000} to
the function
\[F(x_1,\ldots,x_{k\nu})=F_1(x_1,\ldots,x_\nu)+F_2(x_{\nu+1},\ldots,x_{2\nu})+
\cdots+F_k(x_{(k-1)\nu+1},\ldots,x_{k\nu})\]
which is $\b_{i,j}$-periodic in the variable $x_{(i-1)\nu+j}$.

\endproof

The polynomial version of Proposition \ref{th311} is our next result.

\begin{definition} {\rm We say that a system
\[(P_1,P_{1,1},\ldots,P_{1,\nu}),\ldots,(P_d,P_{d,1},\ldots,P_{d,\nu})\]
 of polynomial $(\nu+1)$-tuples  is totally
$\Q$-independent, if for any non-empty subset $\L\subset \{1,\ldots,d\}$
 it is true, that if $a_\l$, $\l\in \L$ are {\it nonzero} integers, then
 the polynomial $\sum_{\l\in \L} a_\l P_\l$ cannot be written as
 a linear combination with rational coefficients of the polynomials
 $P_{\l,j}$, $\l\in \L,1\le j\le\nu$, and of a polynomial with integer coefficients.}
 \end{definition}

\

The most general result on distribution in this paper is the following.
\begin{th} \label{th3110} Let $F_1,\ldots,F_d$ be continuous functions
of $\nu$ variables such that each $F_i$ is periodic in each of its variables with
periods
$\b_{i,1},\ldots,\b_{i,\nu}$, respectively. Let furthermore, $P_1,\ldots,P_d$ and $P_{i,j}$, $1\le i\le d$, $1\le j\le \nu$,
 be polynomials without constant term.
If the system
\[(P_1,P_{1,1}/\b_{1,1},\ldots,P_{1,\nu}/\b_{1,\nu}),\ldots,(P_d,P_{d,1}/\b_{d,1},\ldots,P_{d,\nu}/\b_{d,\nu})\]
 is totally $\Q$-independent, then
the sequence
\be \Bigl(P_1(n)+F_1\bigl(P_{1,1}(n)), \ldots,P_{1,\nu}(n)\bigr),\ldots,
P_d(n)+F_d\bigl(P_{d,1}(n),\ldots,P_{d\nu}(n)\bigr)\Bigr),\label{vvec1}\ee
$n=1,2,\ldots,$ of vectors is uniformly distributed on $[0,1)^d$ {\rm (mod 1)}.
\end{th}

See Proposition \ref{Th43} for the sharpness of this theorem.

\proof. The claim follows from Proposition \ref{th312} exactly as Proposition \ref{th311} resp., \ref{th311*} followed
from Proposition \ref{th31} resp., Proposition \ref{th311}, just use again that uniform distribution of a sequence (\ref{df}) is equivalent to
the uniform distribution of all the sequences (\ref{df1}), and apply  Proposition \ref{th312}.\endproof

\sect{Exactness of the conditions for uniform distribution}
In this section we show that the conditions for uniform distribution set forth in the
preceding section are exact. Let us start by showing that the conditions in Proposition \ref{th3}
and Proposition \ref{th31} cannot be relaxed.
\begin{prop}\label{th4} If $\a$ and $\b_1,\ldots,\b_\nu$ are given and
$\a$ can be written as a
rational linear combination of the numbers
$1,1/\b_1,\ldots,1/\b_\nu$, then there are $\b_j$-periodic continuous functions $f_j$ for $j=1,\ldots,\nu$,
such that the sequence $n\a+f_1(n)+\cdots +f_\nu(n)$, $n=1,2,\ldots$,
is not uniformly distributed {\rm (mod 1)}.
\end{prop}

This proposition does not tell us if, in the case when uniform distribution does not happen,
the sequence has a distribution or not. We shall briefly discuss that question in the next section.

\proof Proposition \ref{th4}. First of all, we may assume that $\a$ is irrational, for
otherwise the claim is trivial by setting $f_j\equiv 0$.

Let $\g_j=1/\b_j$. By the assumption there are integers $q>0$
and $p_0,p_1\ldots,p_d$ such that
\be \a+\sum_{j=1}^\nu \frac{p_j}{q}\g_j=\frac{p_0}{q}.\label{uu}\ee
Multiply this equality by $n$ and write each number on the left as its integral and fractional part
to obtain
\[[n\a]+\{n\a\}+\sum_{j=1}^\nu \frac{p_j}{q}([n\g_j]+\{n\g_j\})=\frac{p_0}{q},\]
from which we can see that
\[\{n\a\}+\sum_{j=1}^\nu \frac{p_j}{q}\{n\g_j\}=\frac{s}{q},\]
where $s$ is an integer. Since on the left we have a linear combination of
numbers lying in $[0,1)$, actually
\[|s|<|q|+|p_1|+\cdots+|p_\nu|=:M.\]

 Therefore,
for every $N$ there is an integer  $t\in [-M,M]$ such that
for at least $N/(2M+1)$ of the numbers $1\le n\le N$ we have
\be \{n\a\}+\sum_{j=1}^\nu \frac{p_j}{q}\{n\g_j\}=\frac{t}{q},\label{hh}\ee
and we may assume that $t$ is the same for infinitely many $N$.

Let $\e>0$ be a small number. Since $\a$ is irrational, by the uniform distribution
of the sequence $\{n\a\}$, $n=1,2,\ldots$, for large $N$ at least $(1-2\e)N$ of the
numbers $\{n\a\}$, $n=1,2,\ldots,N$, belong to $[0,1-\e]$. In a similar vein, if for some
$j$ the number $\g_j$ is irrational, then for large $N$ at least $(1-2\e)N$ of the
numbers $\{n\g_j\}$, $n=1,2,\ldots,N$, belong to $[0,1-\e]$. On the other hand,
if $\g_j$ is rational, then for small $\e>0$ all the
numbers $\{n\g_j\}$, $n=1,2,\ldots,N$, belong to $[0,1-\e]$.
Thus, for at least $(1-2(\nu+1)\e)N$ of the numbers $1\le n\le N$ these relations
hold simultaneously. Choose $N$ sufficiently large so that all these are satisfied.
If, in addition, $2(\nu+1)\e<1/2(2M+1)$, then aut of these
$\ge (1-2(\nu+1)\e)N$ numbers $1\le n\le N$ at least $N/2(2M+1)$ also satisfy (\ref{hh})
(recall that (\ref{hh}) was true for at least $N/(2M+1)$ of the $1\le n\le N$).

Let now $g_j$, $j=1,\ldots,\nu$, be 1-periodic continuous functions such that
\[g_j(x)=\frac{p_j}{q}x\qquad\mbox{for $x \in [0,1-\e]$}.\]
In view of (\ref{hh}), for at least $N/2(2M+1)$ of the numbers $1\le n\le N$
we have
\[n\a+\sum_{i=1}^d g_i(n\g_i)=\frac{t}{q}\qquad \mbox{(mod 1)},\]
and since this is true for infinitely many $N$, it follows that the numbers
on the left for $n=1,2,\ldots$ are not uniformly distributed (mod 1)
(the uniform distribution does not have point masses).

We can complete the proof by setting $f_j(x)=g_j(\g_j x)=g_j(x/\b_j)$\endproof

Next we show that the condition given in Proposition \ref{th312} is sharp.

\begin{prop} \label{th41} If  the numbers $\b_1,\ldots,\b_\nu$  and the polynomials $P_0,P_1,\ldots,P_\nu$
are given and $P_0$ can be written as a rational linear combination of the polynomials
$P_1/\b_1,\ldots,P_\nu/\b_\nu$ and of a polynomial with rational coefficients, then for $j=1,\ldots,\nu$
there are periodic continuous functions $f_j$ of period $\b_j$ such that
 sequence
 \[P_0(n)+f_1(P_1(n))+\cdots +f_\nu(P_\nu(n)), \quad n=1,2,\ldots,\]
  is not
uniformly distributed {\rm (mod 1)}.
\end{prop}
\proof. The proof is almost identical to that of Proposition \ref{th4}.
We may assume that $P_0$ has an irrational coefficient (otherwise just set
$f_j\equiv 0$).

Now start with a relation
\be P_0(x)+\sum_{j=1}^d \frac{p_j}{q}P_j(x)\g_j=\frac{p_0}{q}R_j(x),\label{GG}\ee
where $\g_j=1/\b_j$ and $R_j$ is a polynomial with integer coefficients.
This is the analogue of (\ref{uu}), and from here proceed as we reasoned from (\ref{uu}).
We only need to
mention that since $P_0$ has at least one irrational coefficient, by Weyl's theorem,
for large $N$, at least $(1-2\e)N$ of the
numbers $\{P_0(n)\}$, $n=1,2,\ldots,N$, belong to $[0,1-\e]$, and similarly, if for some
$j$ the polynomial  $P_j\g_j$ has an irrational coefficient,
then for large $N$ at least $(1-2\e)N$ of the
numbers $\{P_j(n)\g_j\}$, $n=1,2,\ldots,N$, belong to $[0,1-\e]$. On the other hand,
if all coefficients of $P_j\g_j$ are rational, then for sufficiently small $\e>0$ all the
numbers $\{P_j(n)\g_j\}$, $n=1,2,\ldots,N$, belong to $[0,1-\e]$. \endproof

The conditions set forth in the vector case in Proposition \ref{th311}
are also best possible:

\begin{prop} \label{Th42} Let $\a_1,\ldots,\a_d$ and
$\b_1,\ldots,\b_d$ be given numbers.
If the pairs $(\a_1,1/\b_1),\ldots,(\a_d,1/\b_d)$ are not totally $\Q$-independent,
 then there are $\b_j$-periodic continuous functions $f_j$ such that
the sequence {\rm (\ref{vvec})} is not uniformly  distributed {\rm (mod 1)} on $[0,1)^d$.
\end{prop}

\proof. We may assume that $1,\a_1,\ldots,\a_d$ are rationally independent (if this is not the
 case, then for $f_j\equiv 0$ the sequence (\ref{vvec}) is not uniformly distributed
 by Weyl's theorem).

 By assumption there is a non-empty set $\L\subset \{1,\ldots,d\}$ and
non-zero integers $a_\l$, $\l\in \L$, such that
\be \sum_{\l\in \L}a_\l\a_\l+\sum_{\l\in \L}\frac{p_\l}{q} \g_\l=\frac{p_0}{q}\label{HH}\ee
where $\g_\l=1/\b_\l$ and $p_\l,q,p_0$ are integers. This is the analogue of (\ref{hh}), and
if we follow the reasoning after (\ref{hh}) (note that we can, since $\sum_\l a_\l \a_\l$ is irrational by the
rational independence of $1,\a_1,\ldots,\a_d$), then we obtain that there are
continuous functions $h_\l$ of period $\b_\l$ such that the sequence
\[\left(\sum_{\l\in \L}a_\l\a_\l\right)n+\sum_{\l\in \L}h_\l(n),\qquad n=1,2,\ldots,\]
is not (mod 1) uniformly distributed. If we set $f_\l =h_\l/a_\l$, then it follows that
the numerical sequence
\[\left(\sum_{\l\in \L}a_\l\a_\l\right)n+\sum_{\l\in \L}a_\l f_\l(n),\qquad n=1,2,\ldots,\]
is not (mod 1) uniformly distributed on $[0,1)$. But then, by Weyl's theorem, the vector
sequence
(\ref{vvec}) (where, say, we set $f_j\equiv 0$ if $j\not\in \L$) is not uniformly
distributed (mod 1) in $[0,1)^d$, either.\endproof

Finally, we show that Theorem \ref{th3110} is sharp. For simplicity we shall formulate
the sharpness only for the case when the functions $F_i$ are of a single variable.

\begin{th} \label{Th43}. Let $\b_1,\ldots,\b_d$ be given numbers
and let $Q_1,\ldots,Q_d$ and $P_1,\ldots,P_d$ be given polynomials without constant term.
If the system of polynomial pairs $(Q_1,P_1/\b_1),\ldots,(Q_d,P_d/\b_d)$ is not totally
$\Q$-independent, then there are continuous periodic functions $f_j$ of period $\b_j$,  $1\le j\le d$,
such that the sequence
\be \Bigl\{\bigl(Q_1(n)+f_1(P_1(n)), \ldots,Q_d(n)+f_d(P_d(n))\bigr)\Bigr\}_{n=1}^\i\label{vvec1**}\ee
is not uniformly distributed {\rm (mod 1)} on $[0,1)^d$.
\end{th}

\proof. We may assume that $Q_1,\ldots,Q_d$ are $\Q$-independent (otherwise
set $f_j\equiv 0$ and apply Weyl's theorem).

By assumption there is a non-empty set $\L\subset \{1,\ldots,d\}$ and
non-zero integers $a_\l$, $\l\in \L$, such that with $\g_\l=1/\b_\l$ we have
\[\sum_{\l\in \L}a_\l Q_\l(x)+\sum_{\l\in \L}\frac{p_\l}{q} P_\l(x)\g_\l=\frac{p_0}{q}R(x),\]
where $p_0, p_\l,q$ are integers and the polynomial $R$ has integer coefficients.
 This is now the analogue of (\ref{HH}) and (\ref{GG}), and from here the
 reasoning is the same that was given after (\ref{HH}) (taking into account
 the necessary modifications in the polynomial case as were given after (\ref{GG})
 and taking also into account that $\sum a_\l Q_\a$ has
 an irrational coefficient since $Q_1,\ldots,Q_d$ are $\Q$-independent):
 by Proposition \ref{th41} there are
continuous functions $h_\l$ of period $\b_\l$ such that the sequence
\[\left(\sum_{\l\in \L}a_\l Q_\l(n)\right)+\sum_{\l\in \L}h_\l(P_\l(n)),\qquad n=1,2,\ldots,\]
is not (mod 1) uniformly distributed. Now set $f_\l =h_\l/a_\l$, and complete the proof
as in Proposition \ref{Th42}.\endproof

\sect{Density and distribution on the torus}\label{secttorus}
In this section we change somewhat our perspective. So far we considered modulo 1 values
as fractional parts on the interval $[0,1)$.  But we can consider them also on
the torus $\R/\Z$, i.e., when we identify the points 0 and 1 in $[0,1]$.
All our results so far
on the (mod 1) density and uniform distribution hold also on the torus without any change. Until now we
were considering sequences (or vector sequences) of the form
\be P_0(0)+F(P_1(n),\ldots,P_\nu(n)),\quad n=1,2,\ldots,\label{tor1}\ee
where, in some cases, $F$ was allowed to be almost periodic. If we only restrict
our attention only to periodic $F$, then it is natural to work with mappings
$G:(\R/\Z)^\nu\to \R/\Z$ (or with
$G:(\R/\Z)^\nu\to (\R/\Z)^d$ in the vector case) and consider
sequences
\be G(P_1(n),\ldots,P_\nu(n)),\quad n=1,2,\ldots\label{tor2}\ee
on $\R/\Z$ (or on $(\R/\Z)^d$ in the vector case), which leads to a somewhat more
general setting.
Indeed, if the $F$ in (\ref{tor1}) is periodic in each variable with periods
$\b_1,\ldots,\b_\nu$, respectively, then by setting
\be G(x_0,x_1,\ldots,x_\nu)=x_0+F(\b_1x_1,\ldots,\b_\nu x_\nu),\quad (x_0,x_1,\ldots,x_\nu)\in (\R/\Z)^{\nu+1},
\label{tor3}\ee
and $P_j^*(x)=P_j(x)/\b_j$, we have
\[P_0(0)+F(P_1(n),\ldots,P_\nu(n))=G(P_0(n),P_1^*(n),\ldots,P_\nu^*(n))\qquad \mbox{(mod 1)},\]
so (\ref{tor2}) include, indeed, (mod 1) sequences of the form (\ref{tor1}).

Any function $G:(\R/\Z)^\nu\to \R/\Z$ can also be considered as a function $G:\R^\nu\to \R/\Z$ by stipulating
$G(x_1,\ldots,x_\nu):=G(\{x_1\},\ldots,\{x_\nu\})$, and this is indeed how we interpreted the values in (\ref{tor2}).
The torus $\R/\Z$ can also be identified with the unit circle under the mapping $x\to e^{2\pi i x}$.
This allows us to speak of the rotation (or winding) number $w$ of a continuous mapping $g:\R/\Z\to \R/\Z$ which
is defined
as the total change of the argument in $e^{2\pi i g(x)}$ as $x$ runs trough the interval $[0,1]$. Intuitively,
$w$ tells us how many times $e^{2\pi i g(x)}$ circles the origin as $e^{2\pi i x}$ makes one full circle
in the positive (counterclockwise) direction. For example, if $m$ is an integer, then the mapping
$x\to mx$ has rotation number $m$, and all $g:\R/\Z\to \R/\Z$ with rotation number $m$
is homotope with the mapping $x\to mx$.

If $G:=(\R/\Z)^\nu\to \R/\Z$ is continuous and we fix $x_2,\ldots,x_\nu$, then $g(x_1)=G(x_1,x_2,\ldots,x_\nu)$ is a mapping
from $\R/\Z$ into itself, and let $w_1$ be its rotation number. Since this is an integer which changes
continuously as we change $x_2,\ldots,x_\nu$ continuously, it follows that $w_1$ is independent of how $x_2,\ldots,x_\nu$ are
fixed, and we call $w_1$  the rotation number of $G$ with respect to $x_1$. The rotation number $w_j$ with respect
the variable $x_j$. $2\le j\le \nu$, is defined similarly. For example, the mapping $G$ in (\ref{tor3})
has rotation number $w_0=1$ with respect to the variable $x_0$ and (because of the periodicity of the function $F$)
rotation number 0 ($w_j=0$) with respect to all other variables $x_j$, $1\le j\le \nu$.

\subsection{Density on the torus}
With the just given definitions we are now ready to extend Proposition \ref{th2pol1} (at least for periodic functions).

\begin{prop} \label{th2pol1tor} Let  $G:=(\R/\Z)^\nu\to \R/\Z$ be a continuous
function with rotation numbers $w_1,\ldots,w_\nu$, respectively, and
let $P_1,\ldots,P_\nu$ be arbitrary polynomials without constant terms. If $w_1P_1+\cdots+w_\nu P_\nu$ has at least
one irrational coefficient, then the sequence
\be \Bigl\{G(P_1(n),\ldots,P_\nu(n))\Bigr\}_{n=1}^\i,\label{001tor}\ee
is dense in $\rz$.
\end{prop}

Note that this is exact in the following sense: if all coefficients of
$w_1P_1+\cdots+w_\nu P_\nu$ are rational, then for
\[G(x_1,\ldots,x_\nu)=w_1x_1+\cdots+w_\nu x_\nu\]
(which clearly has rotation numbers $w_1,\ldots,w_\nu$), the sequence in (\ref{001tor})
is not dense in $\rz$ (for it has only finitely many different terms).

\proof Proposition \ref{th2pol1tor}. The claim easily follows from Proposition \ref{th2pol1}.
In fact, let
\[F(x_1,\ldots,x_\nu)=G(x_1,\ldots,x_\nu)-(w_1x_1+\cdots w_\nu x_\nu).\]
Since all rotation numbers of this $F$ are zero,  $F$ (considered  as a function from $\R^\nu$ to
$\R$) is a 1-periodic function in each of its variable. Furthermore,
if we set
\[P_0:=w_1P_1+\cdots+w_\nu P_\nu,\]
then for all $n$ we have
\[G(P_1(n),\ldots,P_\nu(n))=P_0(n)+F(P_1(n),\ldots,P_\nu(n)),\]
hence the claim in the proposition is a consequence of Proposition
\ref{th2pol1}.\endproof

In a similar vein can one obtain the vector case.

\begin{th}\label{thmaintor} Let $G_i(x_{i,1},\ldots x_{i,\nu})$, $1\le i\le d$, be
continuous functions of $\nu$ variables on $\rz$, and let
$P_{i,j}$, $1\le j\le \nu$, $1\le i\le d$, be polynomials without constant term.
If the rotation number of $G_i$ with respect to the variable
$x_{i,j}$ is $w_{i,j}$ and if the polynomials
\be P_i=w_{i,1}P_{i,1}+\cdots w_{i,\nu}P_{i,\nu},\qquad i=1,\ldots,d,\label{tor4}\ee
are ${\bf Q}$-independent, then the sequence of vectors
\be \Bigl(G_i\bigl(P_{i,1}(n), \ldots,P_{i,\nu}(n)\bigr)\Bigr)_{i=1}^d,\quad n=1,2,\ldots
\label{tor5}\ee
is dense in $(\rz)^d$.
\end{th}
This is exact again, for if the polynomials (\ref{tor4}) are not ${\bf Q}$-independent, i.e.
some linear combination of them with rational coefficients has rational coefficients,
then for
\[G_i(x_{i,1},\ldots,x_{i,\nu})=w_{i,1}x_{i,1}+\cdots w_{i,\nu} x_{i,\nu},\qquad i=1,\ldots,d,\]
the sequence (\ref{tor5}) is not dense in $(\rz)^d$ by Weyl's theorem.

\proof. If we set
\[F_i(x_{i,1},\ldots,x_{i,\nu})=G_i(x_{i,1},\ldots,x_{i,\nu})-(w_{i,1}x_{i,1}+\cdots w_{i,\nu} x_{i,\nu}),\]
then the claim follows from Theorem \ref{thmain} exactly as we deduced Proposition \ref{th2pol1tor}
from Proposition \ref{th2pol1}.\endproof

\subsection{Existence of distribution}
We have seen that when considering the Weyl theory with periodic perturbations, we do not always get
uniform distribution. This raises the question if in those cases the distribution of the sequences
in question exist at all. We shall see that the answer depends on if we work on the torus or on
$[0,1)$ with fractional parts.

Let $X_1,X_2,\ldots$ be a sequence on the torus $({\rz})^d$. We say that this sequence has
distribution $\s$ if for every continuous function $f:(\rz)^d\to \R$ we have
\be \lim_{n\to\i}\frac{1}{n}\sum_{k=1}^n f(X_k)=\int f\;d\s.\label{int1}\ee
An equivalent formulation is that for a dense set\footnote{In fact, for all $a_i,b_i$
for which the boundary of the box $\prod_{i=1}^d[a_i,b_i]$ has zero $\s$-measure.}  of $0< a_i<b_i<1$, $1\le i\le d$,
we have
\be \lim_{n\to\i}\#\left\{1\le k\le n\sep X_k\in \prod_{i=1}^d[a_i,b_i]\right\}=\s\left(\prod_{i=1}[a_i,b_i]\right).\ee

Uniform distribution is when $\s$ is the Haar measure on $(\rz)^d$ (which is the $d$-fold product of
the normalized arc measure on the unit circle). In proving the existence of such a $\s$ it
is sufficient to show that the limit
\[L_f:=\lim_{n\to\i}\frac{1}{n}\sum_{k=1}^n f(X_k)\]
exists for all such $f$, for then $L_f$ is clearly a positive linear functional
and the Riesz representation theorem (cf. \cite[Theorem 2.14]{Rudin}) gives (\ref{int1}) with some measure $\s$.

Now on the torus the existence of the distribution for polynomial sequences holds in all situations.
\begin{th} \label{thtor2} Let $G_1,\ldots,G_d$ be continuous  mappings of $(\rz)^\nu$ into $\rz$,
and let  $P_{i,j}$, $1\le i\le d$, $1\le j\le \nu$,
 be arbitrary polynomials. Then the sequence
\[\Bigl(G_1\bigl(P_{1,1}(n)), \ldots,P_{1,\nu}(n)\bigr),\ldots,
G_d\bigl(P_{d,1}(n),\ldots,P_{d,\nu}(n)\bigr)\Bigr)_{n=1}^\i\]
of vectors has a distribution.
\end{th}
\proof. By introducing $(d-1)\nu$ dummy variables in each $G_i$ (and hence
replacing $\nu$ by $\nu d$) and by listing
$P_{i,j}$ into a single sequence $P_j$ it is sufficient to consider the case when $P_{i,j}$
does not depend on $i$, i.e., it is sufficient to consider sequences
of the form
\[\Bigl(G_1\bigl(P_1(n)), \ldots,P_\nu(n)\bigr),\ldots,
G_d\bigl(P_1(n),\ldots,P_\nu(n)\bigr)\Bigr)_{n=1}^\i.\]

As before, we may assume that each $P_j$ is without a constant term.
Choose a maximal ${\bf Q}$-independent subset of $P_1,\ldots,P_\nu$. We may assume
that this is $P_1,\ldots,P_\mu$ for some $\mu$, and first consider the case when
$\mu\ge 1$. Every $P_j$, $\mu<j\le \nu$,  can be expressed as
\[P_j=R_j+\sum_{s=1}^\mu c_{s,j}P_s,\]
where the $c_{s,j}$ and the coefficients of the polynomial $R_j$ are rational numbers.
Choose an integer $q>0$ such that every $qc_{s,j}$ is an integer, and all coefficients of
all $qR_j$
are integers. It is sufficient to show that the distribution of every subsequence
\[\Bigl(G_1\bigl(P_1(k)), \ldots,P_\nu(k)\bigr),\ldots,
G_d\bigl(P_1(k),\ldots,P_\nu(k)\bigr)\Bigr)_{k=1.\;k\in {\cal N}_r}^\i\]
exists, where ${\cal N}_r=\{k\sep k=ql+r\}_{l=1}^\i$ and $0\le r<q$ is a fixed integer.
If $k\in {\cal N}_r$, then
\[P_j(k)=R_j(r)+\sum_{s=1}^\mu c_{s,j}P_s(k)\qquad \mbox{(mod 1)},\]
so if
\[G_i^*(x_1,\ldots,x_\mu)=G_i\Bigl(x_1,\ldots,x_\mu, R_{\mu+1}(r)+\sum_{s=1}^\mu c_{s,\mu+1}x_s,\ldots,
 R_{\nu}(r)+\sum_{s=1}^\mu c_{s,\nu}x_s\Bigr),\]
 then
\[ G_i(P_1(k)), \ldots,P_\nu(k))=G_i^*(P_1(k)), \ldots,P_\mu(k))\]
for all $i$.

Now if $f$ is any continuous function on $(\rz)^d$ and we set
\[f^*(x_1,\ldots,x_\mu)=f(G_1^*(x_1,\ldots,x_\mu),\ldots,G_d^*(x_1,\ldots,x_\mu)),\]
then we have to show the existence of the limit
\[\lim_{n\to\i}\frac{q}{n}\sum_{k=1,\; k\in {\cal N}_r}^n f^*(P_1(k),\ldots,P_\mu(k)).\]
Since, by Weyl's theorem, $\{(P_1(k),\ldots,P_\mu(k))\}_{k=1,\; k\in {\cal N}_r}^\i$
is uniformly distributed in $(\rz)^\mu$ (a consequence of the ${\bf Q}$-independence of
the polynomials $P_1,\ldots,P_\mu$\footnote{Note that
any ${\bf Q}$-dependence of $P_j(qx +r)$, $1\le j\le \mu$, would automatically extend to a
${\bf Q}$-dependence of $P_1,\ldots,P_\mu$}), the preceding limit exists (and is actually
$\int f^*d\s$ where $\s$ is the Haar measure on $(\rz)^\mu$).

If $\mu=0$, i.e., when all $P_i$ have rational coefficients, then for a $q$ for which all $qP_j$ have integer
coefficients and for $k=ql+r$ we have $P_j(k)=P_j(r)$ (mod 1), i.e.
$(P_1(k),\ldots,P_\nu(k))$ is a constant vector in $(\rz)^\nu$ for $k\in {\cal N}_r$, from which the
claim in the theorem immediately follows.\endproof

Finally, we show that the situation is different if we consider distribution not on the torus $\rz$, but the distribution of the
fractional parts (that belong to $[0,1)$. The following proposition shows that in this case the distribution
may not exist in very simple situations.

\begin{prop}\label{prop} Let $\a$ be irrational. Then there is a  $1/\a$-periodic continuous function $f$
such that the distribution of the fractional parts $\{n\a+f(n)\}$, $n=1,2,\ldots$, does not exist.
\end{prop}

The non-existence of the distribution in this case will mean that for all $b\in [7/8,15/16]$ the limit
\[\lim_{n\to\i}\frac{1}{n}\#\{{1\le k\le n\sep \{k\a+f(k)\}\in [0,b]}\}\]
does not exist.

Recall, however, that if $f$ is $\b$-periodic for some irrational $\b$ and $1,\a$ and $\b$ are rationally independent, then
this phenomenon cannot happen since then the fractional parts $\{n\a+f(n)\}$, $n=1,2,\ldots$, are uniformly distributed
(see Proposition \ref{th3}).

\proof Proposition \ref{prop}. It will be sufficient to construct a 1-periodic function $h$
for which the sequence
$\{n\a+h(n\a)\}$, $n=1,2,\ldots$, of fractional parts has no distribution on $[0,1)$,
for then $f(x):=h(\a x)$ is suitable in Proposition \ref{prop}.

Let $x_m=\{m\a\}$, and $S=\{x_m\}_{m=1}^\i$.
Let $h_0$ be the 1-periodic continuous function
for which
\[h_0(x)=\left\{
  \begin{array}{lll}
    -x & \mbox{$x\in [0,1/2]$}, \\
    -1/2+2(x-1/2) & \mbox{$x\in (1/2,3/4]$},\\
    0&\mbox{$x\in (3/4,1)$}.
  \end{array}\right.\]
Then
\[x+h_0(x)=\left\{
  \begin{array}{lll}
    0 & \mbox{$x\in [0,1/2]$}, \\
    3(x-1/2) & \mbox{$x\in (1/2,3/4]$},\\
    x&\mbox{$x\in (3/4,1)$}.
  \end{array}\right.\]

For $k=0,1,\ldots$ we shall define a 1-periodic continuous function $h_k$ and numbers $M_{k-1}<N_k<M_k$, where
$-1/16<h_{k+1}\le h_k$ will be true for all $k$, and $h_{k+1}$ will be obtained from $h_k$
by modifying $h_k$ on certain finitely many disjoint open intervals $(a_{k,s},b_{k,s})\subset(0,1/2)$, $s\in L_k$,
 over which $h_k(x)=-x$. The modification will be done with the "triangle" function
\[\chi_{k,s}(x)=\left\{
  \begin{array}{lll}
    a_{k,s}-x & \mbox{$x\in (a_{k,s},c_{k,s}]$}, \\
    x-b_{k,s} & \mbox{$x\in (c_{k,s},b_{k,s})$},\\
    0&\mbox{elsewhere in $[0,1)$,}
  \end{array}\right.\]
where $c_{k,s}$ is the midpoint of $(a_{k,s},b_{k,s})$,  i.e. we set
  \be h_{k+1}=h_k+\sum_{s\in L_k} \chi_{k,s}(x),\label{io}\ee
and of course continue this $h_{k+1}$ periodically with period 1.
 We shall also require that
the endpoints $a_{k,s},b_{k,s}$ do not belong to the set $S$, and
$h_{k+1}(x)=h_0(x)$ outside a set consisting of finitely many closed intervals $\subset (0,1/2)$ of total length $<1/8$.
 Note that with such a construction the functions $\{h_k\}_{k=1}^\i$
 uniformly converge monotone
 decreasingly to a continuous function $h$ such that
$x+h(x)=x$ on $[3/4,1)$, and
for $x\in [0,1)$ the fractional part $\{x+h(x)\}$ belongs to the interval
$(7/8,1)$ precisely if $x\in (7/8,1)$ or if $x+h(x)<0$, and in the
latter case the fractional part $\{x+h(x)\}$ actually belongs to $(15/16,1)$.

$h_0$ has already been given, and set $N_0=1$, $M_0=2$. Suppose $h_k,N_k,M_k$ have already
been defined with the property that $h_k(x)=h_0(x)$ outside a set of closed intervals $\subset(0,1/2)$
of total length $<1/8$, say of
total length $1/8-\d_k$.
The measure of the set
\be \{x\in [0,1)\sep \{x+h_k(x)\}\in (7/8,1)\}\label{set}\ee
is at most the measure of $(7/8,1)$ plus the measure of those $x\in [0,1/2]$ where $h_k(x)<h_0(x)$,
hence (\ref{set}) has measure $<1/4$.
Since $h_k$ is 1-periodic, we have
\[\{m\a+h_k(m\a)\}=\{x_m+h_k(x_m)\},\]
and we get from the uniform distribution of the sequence $\{x_m\}_{m=1}^\i$ in $[0,1)$
that there is an $N_{k+1}>M_k$ such that
\be \#\left\{1\le m\le N_{k+1}\sep x_m+h_k(x_m)\in (7/8,1)\ \ \mbox{(mod 1)}\right\}\le \frac{N_{k+1}}{4}.\label{kl}\ee
For the same reason there is an $M_{k+1}>N_{k+1}$ such that the set
\[L_k=\left\{N_{k+1}< m\le M_{k+1}\sep x_m\in (0,1/2), \ x_m+h_k(x_m)=0\right\}\]
has at least
\be \#L_k\ge \left(\frac12-\frac18\right)M_{k+1}\label{13}\ee
elements. For each such $x_m$ we have $h_k(x)=-x$ in a neighborhood of $x_m$ (this follows from
the fact that $a_{l,s},b_{l,s}\not\in S$ for $l<k$ by the induction hypothesis and by the fact
that $h_k(x)<-x$ only on the unions of the intervals $(a_{l,s},b_{l,s})$ with $l<k$), so we can choose
for each $s\in L_k$ an interval $(a_{k,s},b_{k,s})\subset(0,1/2)$ with center at $x_s$ such that
\begin{itemize}
\item $a_{k,s},b_{k,s}\not \in S$,
\item $(a_{k,s},b_{k,s})$ does not contain any of the points $x_m$, $1\le m\le N_{k+1}$,
\item the total length of these $(a_{k,s},b_{k,s})$, $s\in L_k$ is smaller than $\d_k/16(k+1)$,
\item $h_k(x)=-x$ on $(a_{k,s},b_{k,s})$ and
\item these intervals
are pairwise disjoint for different $s$.
\end{itemize}
Finally set $h_{k+1}$
as in (\ref{io}), and this finishes the induction.

The sequence $\{h_k\}_{k=1}^\i$ converges decreasingly and uniformly to a 1-periodic continuous function $h$, and by
the construction we have $h(x_m)=h_k(x_m)$ for all $1\le m\le N_{k+1}$. Hence, (\ref{kl}) yields
\be \liminf_{k\to\i}\frac{1}{N_{k+1}}\#\left\{1\le m\le N_{k+1}\sep x_m+h(x_m)\in (7/8,1)\ \ \mbox{(mod 1)}\right\}\le \frac{1}{4}.\label{kl1}\ee
On the other hand, for all $m\in L_k$ we have
\[\frac{1}{16}\le x_m+h(x_m)\le x_m+h_{k+1}(x_m)<0,\]
 and hence
for all such $m$ the fractional part $\{x_m+h(x_m)\}$ belongs to $(15/16,1)$. But then
(\ref{13}) implies
\be \limsup_{k\to\i}\frac{1}{M_{k+1}}\#\left\{1\le m\le M_{k+1}\sep x_m+h(x_m)\in (15/16,1)\ \ \mbox{(mod 1)}\right\}\ge
\frac12-\frac18.\label{kl2}\ee

Since $(15/16,1)$ is a proper subinterval of $(7/8,1)$, the relations
(\ref{kl1}) and (\ref{kl2}) show that the sequence $\{x_m+h(x_m)\}_{m=1}^\i$,
and hence also the sequence $\Bigl\{\{m\a+h(m\a)\}\Bigr\}_{m=1}^\i$
 (mod 1)  has no distribution.\endproof

\vskip1cm

ELKH-SZTE Analysis and Stochastics Research Group

Bolyai Institute

University of Szeged

Szeged

Aradi v. tere 1, 6720, Hungary
\smallskip

{\it totik@math.u-szeged.hu}

\end{document}